\input amstex
\documentstyle{amsppt}
\magnification 1200
\NoRunningHeads

\topmatter
\title 
A central limit theorem for fields of martingale differences
\endtitle
\author 
Dalibor Voln\'y 
\endauthor
\affil
Laboratoire de Math\'ematiques Rapha\"el Salem, UMR 6085, Universit\'e de Rouen, France
\endaffil

\abstract
We prove a central limit theorem for stationary random fields of martingale differences $f\circ T_{\underline{i}}$, 
$\underline{i}\in \Bbb Z^d$, where $T_{\underline{i}}$ is a $\Bbb Z^d$ action and the martingale is given by a commuting filtration.
The result has been known for Bernoulli random fields;
here only ergodicity of one of commuting transformations generating the $\Bbb Z^d$ action is supposed.
\endabstract

\endtopmatter

\document

\heading{Introduction}
\endheading

In study of the central limit theorem for dependent random variables, the case of martingale difference sequences has played
an important role, cf\. Hall and Heyde, \cite{HaHe}. Limit theorems for random fields of martingale differences were studied 
for example by Basu and Dorea \cite{BD}, Morkvenas \cite{M}, Nahapetian \cite{N}, Poghosyan and Roelly \cite{PR}, Wang and Woodroofe \cite{WaW}.
Limit theorems for martingale differences enable a research of much more complicated processes and random fields. The method
of martingale approximations, often called Gordin's method, originated by Gordin's 1969 paper \cite{G1}. The approximation is 
possible for random fields as well, for most recent results cf\. e.g\. \cite{WaW} and \cite{VWa}. Remark that another approach
was introduced by Dedecker in \cite{D} (and is being used since); it applies both to sequences and to random fields.

For random fields, the martingale structure can be introduced in several different ways.
Here we will deal with a stationary random field $f\circ T_{\underline{i}}$, $\underline{i}\in \Bbb Z^d$, where $f$ is a measurable 
function on a probability space $(\Omega, \mu, \Cal A)$ and $T_{\underline{i}}$, $\underline{i}\in \Bbb Z^d$, is a group of
commuting probability preserving transformations of $(\Omega, \mu, \Cal A)$ (a $\Bbb Z^d$ action). 
By $e_i \in \Bbb Z^d$ we denote the vector $(0,...,1, ...,0)$ having 1 on the $i$-th place and 0 at all other places, $1\leq i\leq d$.\newline
$\Cal F_{\underline{i}}$, $\underline{i} = (i_1\dots,i_d) \in \Bbb Z^d$, is an invariant commuting filtration (cf\. D\. Khosnevisan, 
\cite{K}) if 
\roster
\item"(i)" $\Cal F_{\underline{i}}= T^{-\underline{i}} \Cal F_{\underline{0}}$ for all $\underline{i} \in \Bbb Z^d$,
\item"(ii)" $\Cal F_{\underline{i}} \subset F_{\underline{j}}$ for $\underline{i} \leq \underline{j}$ in the lexicographic order, and
\item"(iii)" $\Cal F_{\underline{i}} \cap \Cal F_{\underline{j}} = \Cal F_{\underline{i}\wedge 
\underline{j}}$, $\underline{i}, \underline{j}\in \Bbb Z^d$, and $\underline{i}\wedge \underline{j} =
(\min\{i_1, j_1\}, \dots, \min\{i_d, j_d\})$.
\endroster
If, moreover, $E\Big( E(f | \Cal F_{\underline{i}}) \big| \Cal F_{\underline{j}}\Big) = E(f | \Cal F_{\underline{i}\wedge \underline{j}})$,
for every integrable function $f$, we say that the filtration is {\it completely commuting} (cf\. \cite{G2}, \cite{VWa}). \newline
By $\Cal F_l^{(q)}$, $1\leq q\leq d$, $l\in \Bbb Z$, we denote the $\sigma$-algebra generated by the union of all $\Cal F_{\underline{i}}$ 
with $i_q\leq l$.
For $d=2$ we by $\Cal F_{\infty, j} = \Cal F_j^{(2)}$ denote the $\sigma$-algebra generated by the union of all $\Cal F_{i, j}$, 
$i\in \Bbb Z$, and in the same way we define $\Cal F_{i, \infty}$.

We sometimes denote $f\circ T_{\underline{i}}$ by $U_{\underline{i}}f$; $f$ will always be from $\Cal L^2$.

We say that $U_{\underline{i}}f$, $\underline{i} \in \Bbb Z^d$, is a {\it field of martingale differences} if $f$ is 
$\Cal F_{\underline{0}}$-measurable and whenever $\underline{i} = (i_1\dots,i_d) \in \Bbb Z^d$
is such that $i_q \leq 0$ for all $1\leq q\leq d$ and at least one inequality is strict then $E(f\,|\,\Cal F_{\underline{i}}) =0$.

Notice that $U_{\underline{i}}f$ is then $\Cal F_{\underline{i}}$-measurable, $\underline{i} = (i_1\dots,i_d) \in \Bbb Z^d$, and if 
$\underline{j} = (j_1\dots,j_d) \in \Bbb Z^d$ is such that $j_k \leq i_k$ for all $1\leq k\leq n$ and at least one inequality is strict, 
$E(U_{\underline{i}} f\,|\,\Cal F_{\underline{j}}) =0$. \newline
Notice that by commutativity, if $U_{\underline{i}}f$ are martingale differences then
$$
  E(f \,|\, \Cal F_{-1}^{(q)}) = 0
$$
for all $1\leq q \leq d$. $(f\circ T_{e_q}^j)_j$ is thus a sequence of martingale differences for the filtration of $\Cal F_{j}^{(q)}$.
In particular, for $d=2$, $(f\circ T_{e_2}^j)$ is a sequence of martingale differences for the filtration of  
$\Cal F_{\infty, j}  = \Cal F_j^{(2)}$.

Recall that a measure preserving transformation $T$ of $(\Omega, \mu, \Cal A)$ is said to be {\it ergodic} if for any $A\in \Cal A$ such 
that $T^{-1}A = A$, $\mu(A) =0$ or $\mu(A)=1$. Similarly, a $\Bbb Z^d$ action $(T_{\underline{i}})_{\underline{i}}$ is ergodic if for any 
$A\in \Cal A$ such that $T_{-\underline{i}}A = A$, $\mu(A) =0$ or $\mu(A)=1$.

A classical result by Billinsley and Ibragimov says that if $(f\circ T^i)_i$ is an ergodic sequence of martingale differences, the central
limit theorem holds. The result does not hold for random fields, however. 

\subheading{Example} As noticed in paper by Wang, Woodroofe \cite{WaW}, for a 2-dimensional 
random field $Z_{i,j} = X_iY_j$ where $X_i$ and $Y_j$, $i,j \in \Bbb Z$, are mutually independent $\Cal N(0,1)$ random variables,
we get a convergence towards a non normal law. The random field of $Z_{i,j}$ can be represented by a non ergodic action of $\Bbb Z^2$:

Let $(\Omega, \mu,\Cal A)$ be a product of probability spaces $(\Omega', \mu',\Cal A')$ and $(\Omega'', \mu'',\Cal A'')$ equipped with 
ergodic measure preserving transformations $T'$ and $T''$. On $\Omega$ we then define a measure preserving $\Bbb Z^2$ action
$T_{i,j}(x,y) = ({T'}^ix, {T''}^jy)$. The $\sigma$-algebras $\Cal A', \Cal A''$ are generated by $\Cal N(0,1)$ iid sequences of random variables 
$(e'\circ {T'}^i)_i$ and $(e''\circ {T''}^i)_i$ respectively. The dynamical systems $(\Omega', \mu',\Cal A', T')$ and 
$(\Omega'', \mu'',\Cal A'', T'')$ are then Bernoulli hence ergodic (cf\. \cite{CSF}). On the other hand, for any 
$A'\in\Cal A'$, $A'\times \Omega''$ is $T_{0,1}$-invariant hence $T_{0,1}$ is not an ergodic transformation. Similarly we get that
$T_{1,0}$ is not an ergodic transformation either. By ergodicity of $T', T''$, $A'\times \Omega''$, $A'\in\Cal A'$, are the only $T_{0,1}$-invariant 
measurable subsets of $\Omega$ and $A''\times \Omega'$, $A''\in\Cal A''$, are the only $T_{1,0}$-invariant measurable subsets of $\Omega$
(modulo measure $\mu$). Therefore, the only measurable subsets of $\Omega$ which are invariant both for $T_{0,1}$ and for $T_{1,0}$
are of measure 0 or of measure 1, i.e\. the $\Bbb Z^2$ action $T_{i,j}$ is ergodic.\newline
On $\Omega$ we define random variables $X, Y$ by $X(x,y) = e'(x)$ and $Y(x,y) = e''(y)$. The random field of $(XY)\circ T_{i,j}$ then
has the same distribution as the random field of $Z_{i,j} = X_iY_j$ described above. The natural filtration of 
$\Cal F_{i,j} = \sigma\{(XY)\circ T_{i',j'} : i'\leq i, j'\leq j\}$ is commuting and $((XY)\circ T_{i,j})_{i,j}$ is a field of martingale
differences.
\medskip

A very important particular case of a $\Bbb Z^d$ action is the case when the $\sigma$-algebra $\Cal A$ is generated by iid random variables
$U_{\underline{i}}e$, $\underline{i} \in \Bbb Z^d$. The $\sigma$-algebras $\Cal F_{\underline{j}} = \sigma\{U_{\underline{i}} : i_k \leq j_k, 
k=1,\dots,d\}$ are then a completely commuting filtration and if  $U_{\underline{i}}f$, $\underline{i} \in \Bbb Z^d$ is
a martingale difference random field, the central limit theorem takes place (cf\. \cite{WW}). This fact enabled to prove 
a variety of limit theorems by martingale approximations (cf\. e.g\. \cite{WaW}, \cite{VWa}).\newline
For Bernoulli random fields, other methods of proving limit theorems have been used, cf\. e.g\. \cite{ElM-V-Wu}, \cite{Wa}, \cite{BiDu}.
\medskip

The aim of this paper is to show that for a martingale difference random field, the CLT can hold under assumptions weaker than
Bernoullicity.
\bigskip

\heading{Main result} 
\endheading

Let $T_{\underline{i}}$, $\underline{i}\in \Bbb Z^d$, be a $\Bbb Z^d$ action of measure preserving transformations
on $(\Omega, \Cal A, \mu)$, $(\Cal F_{\underline{i}})_{\underline{i}}$, $\underline{i}\in \Bbb Z^d$, be a commuting filtration.
By $e_i \in \Bbb Z^d$ we denote the vector $(0,...,1, ...,0)$ having 1 on the $i$-th place and 0 at all other places, $1\leq i\leq d$.


\proclaim{Theorem} Let $f\in L^2$, be such that $(f\circ T_{\underline{i}})_{\underline{i}}$ is a field of martingale differences for a
completely commuting filtration $\Cal F_{\underline{i}}$. If at least one of the transformations $T_{e_i}$, $1\leq i\leq d$, is ergodic 
then the central limit theorem holds, i.e\. for $n_1,\dots, n_d \to\infty$ the distributions of
$$
  \frac1{\sqrt{n_1\dots n_d}} \sum_{i_1=1}^{n_1} \dots \sum_{i_d=1}^{n_d} f\circ T_{(i_1,\dots,i_d)}
$$
weakly converge to $\Cal N(0, \sigma^2)$ where $\sigma^2 = \|f\|_2^2$.
\endproclaim

\noindent  \underbar{\it Remark 1.} The results from \cite{VoWa} remain valid for $\Bbb Z^d$ actions
satisfying the assumptions of the Theorem, Bernoullicity thus can be replaced by ergodicity of one of the 
transformations $T_{e_i}$. Under the assumptions of the Theorem we thus also get a weak invariance principle.
\cite{VoWa} implies many earlier results, cf\. references therin and in \cite{WaW}.

\demo{Proof}

We prove the theorem for $d=2$. Proof of the general case is similar. \newline
We suppose that the transformation $T_{0,1}$ is ergodic and $\|f\|_2=1$. To prove the central limit theorem for the random field it is 
sufficient to prove that for $m_k,n_k \to \infty$ as $k\to\infty$,
$$
  \frac1{\sqrt{m_kn_k}} \sum_{i=1}^{m_k} \sum_{j=1}^{n_k} f\circ T_{i,j}\,\,\,\,\text{converge in distribution to}\,\,\,\, \Cal N(0, 1). \tag1 
$$ 
\medskip

Recall the central limit theorem by D.L\. McLeish (cf\. \cite{M}) saying that if $X_{n,i}$, $i=1,\dots,k_n$, is an array of martingale 
differences such that
\roster
\item"(i)" $\max_{1\leq i\leq k_n} |X_{n,i}| \to 0$ in probability,
\item"(ii)" there is an $L<\infty$ such that $\max_{1\leq i\leq k_n} X_{n,i}^2 \leq L$ for all $n$, and
\item"(iii)" $\sum_{i=1}^{k_n} X_{n,i}^2 \to 1$ in probability, 
\endroster
then $\sum_{i=1}^{k_n} X_{n,i}$ converge to $\Cal N(0, 1)$ in law.
\medskip

Next, we will suppose $k_n=n$; we will denote $U_{i,j} f = f\circ T_{i,j}$. For a given positive integer $v$ and positive integers $u, n$
define
$$
  F_{i,v} = \frac1{\sqrt v} \sum_{j=1}^v U_{i,j} f,\quad X_{n,i} = \frac1{\sqrt n} F_{i,v},\,\,\,\, i=1,\dots,n.
$$ 
Clearly, $X_{n,i}$ are martingale differences for the filtration $(\Cal F_{i, \infty})_i$. We will verify the assumptions of McLeish's theorem.

The conditions (i) and (ii) are well known to follow from stationarity. For reader's convenience we recall their proofs.

(i) For $\epsilon >0$ and any integer $v\geq 1$,
$$\multline
  \mu(\max_{1\leq i\leq n} |X_{n,i}| > \epsilon) \leq \sum_{i=1}^n \mu(|X_{n,i}| > \epsilon) =
  n \mu\Big(\Big| \frac1{\sqrt{nv}} \sum_{j=1}^v U_{0,j} f \Big| > \epsilon \Big) \leq \\
  \leq \frac1{\epsilon^2} E\Big( \Big(\frac1{\sqrt{v}} \sum_{j=1}^M U_{0,j} f\Big)^2 1_{|\sum_{j=1}^v U_{0,j} f|\geq 
  \epsilon \sqrt{nv}}\Big) \to 0
  \endmultline
$$
as $n\to \infty$; this proves (i). Notice that that the convergence is uniform for all $v$.

To see (ii) we note
$$
  \Big(\max_{1\leq i\leq n} |X_{n,i}|\Big)^2 \leq \sum_{i=1}^n X_{n,i}^2 = 
  \frac1n \sum_{i=1}^n \Big(\frac1{\sqrt v} \sum_{j=1}^v U_{i,j} f\Big)^2
$$
which implies $E\Big(\max_{1\leq i\leq n} |X_{n,i}|\Big)^2 \leq 1$.
\medskip\cite{WaW}

It remains to prove (iii). \newline
Let us fix a positive integer $m$ and for constants $a_1,\dots, a_m$ consider the sums
$$
  \sum_{i=1}^m a_i \sum_{j=1}^v U_{i,j} f,\,\,\,\, v\to\infty.
$$

Then $(\sum_{i=1}^m a_i U_{i,j} f)_j$, $j=1,2,\dots$, are martingale differences for the filtration $(\Cal F_{\infty, j})_j$ and by 
the central limit theorem of Billingsley and Ibragimov \cite{Bil}, \cite{I} (we can also prove using the McLeish's theorem)
$$
  \frac1{\sqrt v} \sum_{j=1}^v \big(\sum_{i=1}^m a_i U_{i,j} f\big)
$$ 
converge in law to $\Cal N(0,\sum_{i=1}^m a_i^2)$. Notice that that here we use the assumption of ergodicity of $T_{0,1}$. \newline
From this it follows that the random vectors $(F_{1,v},\dots, F_{m,v})$ where
$$
  F_{u,v} = \frac1{\sqrt v} \sum_{j=1}^v U_{u,j} f, \,\,\,\,u=1,\dots, m,
$$  
converge in law to a vector $(W_1, \dots, W_m)$ of $m$ mutually independent and $\Cal N(0, 1)$ distributed random variables. 
For a given $\epsilon>0$, if $m=m(\epsilon)$ is sufficiently big then we have $\big\|1- (1/m) \sum_{u=1}^m F_{u,v}^2\big\|_1 < \epsilon/2$.
Using a truncation anrgument we can from the convergence in law of $(F_{u,v}, \dots, F_{m,v})$ towards $(W_1, \dots, W_m)$
deduce that for
$m=m(\epsilon)$ sufficiently big and $v$ bigger than some $v(m,\epsilon)$,
$$
  \Big\|1- \frac1m \sum_{u=1}^m F_{u,v}^2\Big\|_1 < \epsilon.
$$
Any positive integer $N$ can be expressed as $N= pm +q$ where $0\leq q\leq m-1$. 
Therefore 
$$
  1 - \frac1{N} \sum_{i=1}^N F_{i,v}^2 = \frac{m}N \sum_{k=0}^{p-1} \Big( \frac1m \sum_{i=km+1}^{(k+1)m} F_{i,v}^2 - 1\Big)
  + \frac1N \sum_{i=mp+1}^N F_{i,v}^2 - \frac{q}N. \tag2
$$  
There exists an $N_\epsilon$ such that for $N\geq N_\epsilon$ we have $\|\frac1N \sum_{i=mp+1}^N F_{i,v}^2 - \frac{q}N \|_1 <\epsilon$ hence 
if $v\geq v(m,\epsilon)$ and $N\geq N_\epsilon$ then
$$
  \Big\|1 - \frac1{N} \sum_{i=1}^N F_{i,v}^2 \Big\|_1 =
  \Big\|1 - \frac1{Nv} \sum_{i=1}^N \big( \sum_{j=1}^v U_{i,j} f \big)^2 \Big\|_1  < 2\epsilon. \tag3
$$  
This proves that for $\epsilon>0$ there are positive integers $v(m,\epsilon/2)$ and $N_\epsilon$ such that for $M\geq v(m,\epsilon/2)$
and $n\geq N_\epsilon$, for $X_{n,i} = (1/\sqrt n) F_{i,M}$
$$
  \big\|\sum_{i=1}^n X_{n,i}^2 -1\big\|_1 = \Big\| \sum_{i=1}^n \big(\frac1{\sqrt{nM}} \sum_{j=1}^M U_{i,j} f\big)^2 - 1\Big\|_1 < \epsilon.
$$  
\medskip

In the general case we can suppose that $T_{e_d}$ is ergodic (we can permute the coordinates). Instead of $T_{i,j}$ we will consider
transformations $T_{\underline{i},j}$ where $\underline{i} \in \Bbb Z^{d-1}$ and in (3), instead of segments $\{km+1, \dots, km+m\}$
we take boxes of $(k_1m+i_1, \dots, k_{d-1}m +i_{d-1})$, $i_1, \dots, i_{d-1} \in \{1,\dots, m\}$.
\medskip

This finishes the proof of the Theorem.

\enddemo
\qed

\bigskip

\noindent  \underbar{\it Remark 2.} For any positive integer $d$ there exists a random field of martingale differences 
$(f\circ T_{\underline{i}})$ for a commuting filtration of $\Cal F_{\underline{i}}$ where $T_{\underline{i}}$, $\underline{i}\in \Bbb Z^d$,
is a non Bernoulli $\Bbb Z^d$ action and all $T_{e_i}$, $1\leq i\leq d$, are ergodic.

To show this we take a Bernoulli $\Bbb Z^d$ action $T_{\underline{i}}$, $\underline{i}\in \Bbb Z^d$ on $(\Omega, \Cal A, \mu)$ generated 
by iid random variables $(e\circ T_{\underline{i}})$ as defined e.g\. in \cite{WaW} or \cite{VWa}.\newline
Then we take another $\Bbb Z^d$ action of irrational rotations on the unit circle (identified with the interval $[0, 1)$) generated by 
$\tau_{e_i} = \tau_{\theta_i}$, $\tau_{\theta_i}x = x+\theta_i$ mod 1; $\theta_i$, $1\leq i\leq d$, are linearly independent irrational numbers.
The unit circle is equipped with the Borel $\sigma$-algebra $\Cal B$ and the (probability) Lebesgue measure $\lambda$.\newline
On the product $\Omega \times [0, 1)$ with the product $\sigma$-algebra and the product measure we define the product $\Bbb Z^d$ action
$(T_{\underline{i}}\times \tau_{\underline{i}})(x,y) = (T_{\underline{i}}x, \tau_{\underline{i}}y)$. Because the product of ergodic transformations
is ergodic, for every $e_i$, $1\leq i\leq d$, $T_{e_i}\times \tau_{e_i}$ is ergodic. The product $\Bbb Z^d$ action is not Bernoulli
(it has irrational rotations for factors).

On $\Omega \times [0, 1)$ we define a filtration
$\Cal F_{(i_1,\dots,i_d)} = \sigma\{U_{(i',\dots,i'_d)}e \circ \pi_1, i'-1\leq i_1, \dots, i'_d\leq i_d, \pi_2^{-1}\Cal B\}$
where $\pi_1, \pi_2$ are the coordinate projection of $\Omega \times [0, 1)$.\newline
The filtration defined above is commuting and we can find a random field of martingale differences satisfying the assumptions of the 
Theorem. 
\medskip

\noindent  \underbar{\it Remark 3.} In the one dimensional central limit theorem, non ergodicity implies a convergence towards a mixture
of normal laws. This comes from the fact that using a decomposition of the measure $\mu$ into ergodic components, we get the ``ergodic case'' 
for each of the components (cf\. \cite{V}); the variance is given by the limit of $(1/n) \sum_{i=1}^n U^if^2$ which by the Birkhoff
Ergodic Theorem exists a.s\. and in $L^1$ and is $T$-invariant. In the case of a $\Bbb Z^2$ action (taking $d=2$ for simplicity), 
the limit for $T_{0,1}$ need not be $T_{1,0}$-invariant. This is exactly the case described in the Example and eventually we got there
a convergence towards a law which is not normal.
\bigskip

\noindent {\bf Acknowledgement.}  I am very thankful to J\'er\^ome Dedecker for his remarks, comments, and encouragement. I am also thankful
to Zemer Kosloff; the idea/conjecture that it is the ergodicity of coordinate factors of the $\Bbb Z^d$ action which can imply the central limit 
theorem came out first in our discussion after my lecture in April 2014.

\enddocument

\end